\documentclass[a4paper]{amsart}

\usepackage{amssymb,latexsym}
\usepackage{amsmath}
\usepackage{amsfonts}
\usepackage{amsthm}
\usepackage{amssymb}

\theoremstyle{plain}

\theoremstyle{definition}

\begin{document}

\title[On elliptic curves induced by rational Diophantine quadruple]{On elliptic curves induced by rational Diophantine quadruples}

\author[A. Dujella]{Andrej Dujella}
\address[Andrej Dujella]{
Department of Mathematics\\
Faculty of Science\\
University of \mbox{Zagreb}\\
Bijeni{\v c}ka cesta 30, 10000 Zagreb, Croatia
}
\email{duje@math.hr}

\author[G. Soydan]{G\"okhan Soydan}
\address[G\"okhan Soydan]{
Department of Mathematics \\
Bursa Uluda\u{g} University \\
16059 Bursa, Turkey
}
\email{gsoydan@uludag.edu.tr}



\thanks{A. D. was supported by the Croatian Science Foundation under the project no. IP-2018-01-1313 and
the QuantiXLie Center of Excellence, a project co-financed by the
Croatian Government and European Union through the European Regional Development Fund - the Competitiveness and Cohesion Operational Programme (Grant KK.01.1.1.01.0004). G. S. was supported by the Research Fund of Bursa Uluda\u{g} University
under Project No: F-2020/8.}
\begin{abstract}
In this paper, we consider elliptic curves induced by rational Diophantine quadruples,
i.e. sets of four nonzero rationals such that the product of any two of them plus $1$ is
a perfect square. We show that for each of the groups
$\mathbb{Z}/2\mathbb{Z} \times \mathbb{Z}/k\mathbb{Z}$ for $k=2,4,6,8$,
there are infinitely many rational Diophantine quadruples with the property that the
induced elliptic curve has this torsion group. We also construct curves with moderately large rank
in each of these four cases.
 \end{abstract}

\maketitle

\section{Introduction}
A rational Diophantine $m$-tuple is a set of $m$ distinct nonzero rationals
with the property that the product of any two of its distinct elements plus $1$ is
a perfect square. The first rational Diophantine quadruple was found by Diophantus himself,
and it was the set $\{\frac{1}{16}, \frac{33}{16}, \frac{17}{4}, \frac{105}{16} \}$. The first Diophantine in integers, the set
$\{1,3,8,120\}$ was found by Fermat, and Euler was able to extended Fermat's set to a
rational Diophantine quintuple by the fifth element $\frac{777480}{8288641}$.
Recently, Stoll \cite{Stoll} proved that Euler's extension of
Fermat's quadruple to the rational quintuple is unique.
In 1999, Gibbs (see \cite{Gibbs2006}) found the first rational sextuple
$\{\frac{11}{192}, \frac{35}{192}, \frac{155}{27}, \frac{512}{27}, \frac{1235}{48}, \frac{180873}{16}\}$,
while in 2017, Dujella, Kazalicki, Miki\'c and Szikszai \cite{DKMS} proved that there
exist infinitely many rational Diophantine sextuples (see also \cite{DK-tichy,DKP-jnt,DKP-split}).
For surveys of results and conjectures concerning Diophantine $m$-tuples and their generalizations see
\cite{D-notices} and \cite[Sections 14.6, 16.7]{D-book}.

Let $\{a, b, c\}$ be a rational Diophantine triple.
In order to extend this triple to a
rational Diophantine quadruple, we need to find a rational number $x$ such
that $ax+1$, $bx+1$ and $cx+1$ are squares of rational numbers. By multiplying
these three conditions, we obtain a single condition
$y^2 = (ax + 1)(bx + 1)(cx + 1)$, which is the equation of an elliptic curve.
We say that this elliptic curve is induced by the rational Diophantine triple $\{a,b,c\}$.
The question of possible Mordell-Weil groups of such elliptic curve over $\mathbb{Q}$, $\mathbb{Q}(t)$ and
quadratic fields, was considered in several papers
(see \cite{ADP,D-rocky,D-glas,DJS-racsam,DM2,DP-lms,DP-racsam,DP-glas,DP-jga,YUB}).
In particular, it is shown in \cite{D-glas} that all four torsion groups that are allowed by Mazur's theorem
for elliptic curves with full $2$-torsion, i.e.
$\mathbb{Z}/2\mathbb{Z} \times \mathbb{Z}/k\mathbb{Z}$ for $k=2,4,6,8$, can appear as a torsion group
of an elliptic curve induced by a rational Diophantine triple (note however that this question is still open
for elliptic curves induced by integer Diophantine triples, see \cite{D-bordo,DM1}).
Dujella and Peral used elliptic curves induced by Diophantine triples to
construct elliptic curves over $\mathbb{Q}(t)$ and $\mathbb{Q}$ with given torsion group and large
rank. Let us mention that the curves with the largest known rank over $\mathbb{Q}$ and torsion groups
$\mathbb{Z}/2\mathbb{Z} \times \mathbb{Z}/k\mathbb{Z}$ for $k=4$, $6$ and $8$ can be induced by rational
Diophantine triples (see e.g. \cite{DP-glas}). The details about the current rank
records and the corresponding curves can be found at the web page \cite{D-torsion}.

In this paper, we will consider similar questions but for elliptic curves induced by rational
Diophantine quadruples. Let $\{a,b,c,d\}$ be a rational Diophantine quadruple.
To extend it to a rational Diophantine quintuples, we need to find a rational number $X$ such
that $aX+1$, $bX+1$, $cX+1$ and $dX+1$ are squares of rational numbers.
As in the case of extension of triple to quadruples, we can multiply
these four conditions. We obtain the equation
$$ Y^2 = (aX + 1)(bX + 1)(cX + 1)(dX+1), $$
which is the equation of a genus $1$ curve.
By the substitution
\[ y=\frac{Y(d-a)(d-b)(d-c)}{(dX+1)^2}\,, \qquad
x=\frac{(aX+1)(d-b)(d-c)}{dX+1}\,, \]
we obtain the following elliptic curve
\begin{equation} \label{4}
E: \quad y^2=x(x+(b-a)(d-c))(x+(c-a)(d-b)).
\end{equation}
We say that this elliptic curve is induced by the rational Diophantine quadruple $\{a,b,c,d\}$.
Such curves were considered in \cite{D-irr}, where several examples of curves with torsion group
$\mathbb{Z}/2\mathbb{Z} \times \mathbb{Z}/2\mathbb{Z}$ and rank equal to $8$ were obtained.

In this paper, we consider the question which torsion groups are possible for elliptic curves
induced by rational Diophantine quadruples. We will show that, similarly as in the cases of
elliptic curves induced by rational Diophantine triples,
all four torsion groups that are allowed by Mazur's theorem,
$\mathbb{Z}/2\mathbb{Z} \times \mathbb{Z}/k\mathbb{Z}$ for $k=2,4,6,8$, are possible,
and in fact they can be achieved for infinitely many rational Diophantine quadruples.
This is the main result of our paper.
In each of four cases, we will also find curves with moderately large rank.

There are three non-trivial rational $2$-torsion points on $E$:
\[ A=(0,0),\quad B=(-(b-a)(d-c),0),\quad C=(-(c-a)(d-b),0), \]
and another two obvious rational points:
\begin{align*}
P &= ((b-a)(c-a),\,(b-a)(c-a)(d-a)), \\
Q &=((ad+1)(bc+1),\, \sqrt{(ab+1)(ac+1)(ad+1)(bc+1)(bd+1)(cd+1)}).
\end{align*}
The points $P$ and $Q$ will play an important role in our constructions.
In particular, we will be interested in question under which assumptions
these points may have finite order.
%
%
%
%

\section{Torsion $\mathbb{Z}/2\mathbb{Z} \times \mathbb{Z}/2\mathbb{Z}$} \label{sec:z2z2}

\subsection{Subquadruples of Gibbs' examples of rational Diophantine sextuples}

In \cite{D-irr}, the first known rational Diophantine sextuples found by Gibbs \cite{Gibbs2006}
were analyzed, and by considering ranks of elliptic curves induced by subquadruples
of these sextuples, several elliptic curves with rank $8$ were found.

In \cite{Gibbs2016}, Gibbs presented over $1000$ examples of rational Diophantine sextuples
of low height which are obtained by a computational search. By a systematic computation of the rank
of elliptic curves induced by subquadruples of these new sextuples, we found several curves with rank $9$.
These curves correspond to the following quadruples:

$$
\{a,b,c,d\} = \left\{\frac{1218560}{31752611}, \frac{111}{77}, \frac{34191}{19712}, \frac{1155}{16}\right\} $$

$$
\{a,b,c,d\} = \left\{\frac{8064}{597529}, \frac{1408}{75}, \frac{16225}{48}, \frac{3337875}{16}\right\}
$$

$$
\{a,b,c,d\} = \left\{\frac{43875}{232324}, \frac{71200}{47961}, \frac{1633824}{623045}, \frac{671}{20} \right\}
$$

\subsection{Parametric families of rational Diophantine quadruples}

To improve the results from the previous subsection, we considered several parametric
families of rational Diophantine quadruples.
The most successful was the following approach, which was used previously
in \cite{D-mixed} for construction of rational Diophantine sextuples with mixed signs.

Let $\{a,b,c\}$ be a regular rational Diophantine triple, i.e. $c=a+b+2r$, where $ab+1=r^2$.
Consider now other regular triple $\{b,c,d\}$, where $d=b+c+2\sqrt{bc+1} = a+4b+4r$.
The only remaining condition in order that $\{a,b,c,d\}$ be a Diophantine quadruple
is that $ad+1$ is a perfect square. This condition leads to $(a+2r)^2-3 = \Box$,
and it is satisfied if we take $r=\frac{u-a}{2}$, where $u=\frac{t^2+3}{2t}$
for $t \in \mathbb{Q}$, $t\neq 0$.
Thus we obtain a two-parametric family of rational Diophantine quadruples $\{a,b,c,d\}$, where
\begin{align*}
b &= \frac{(t^2-2at-4t+3)(t^2-2at+4t+3)}{16 t^2 a}, \\
c &= \frac{(t^2+2at+4t+3)(t^2+2at-4t+3)}{16 t^2 a}, \\
d &= \frac{(t-1)(t+3)(t-3)(t+1)}{4 t^2a},
\end{align*}
and the corresponding family of elliptic curves, which is equivalent to
$$ y^2 = x^3 + A_1x^2 + B_1x, $$
where
\begin{align*}
A_1 &= 6t^8-48t^6a^2+96t^4a^4-120t^6+992t^4a^2+708t^4-432t^2a^2-1080t^2+486, \\
B_1 &= (t^2+2at-1)(t^2-6at-9)(3t^2+2at-3)(t^2-2at-9) \\
& \times (t^2+6at-9)(t^2-2at-1)(t^2+2at-9)(3t^2-2at-3).
\end{align*}
We were able to find several curves with rank $9$ within this family.

However, to find curves with slightly larger rank, we noted that if $a^2+1$ is a perfect square,
then four quadratic factors appearing in the factorization of $B_1$ have square discriminant
with the respect to $t$, and thus factorize into two linear factors. This increases the trivial
upper bound for the rank, $\omega(B_1)+\omega(A_1^2-4B_1)$ (see e.g. \cite[Section 15.5]{D-book}),
and indicates that the rank indeed could be larger.
Thus we take
$$ a= \frac{v^2-1}{2v}, $$
and search for high-rank curves with reasonably small numerators and denominators
of the parameters $t$ and $v$. We used the standard sieving methods based on Mestre-Nagao sums
(see \cite{Mestre,Nagao}) and computing Selmer rank as an upper bound for the rank.
The Selmer ranks and exact ranks (for best candidates for high rank) are computed by {\tt mwrank} \cite{Cremona}.
We found two curves with rank $10$, corresponding to
$t=\frac{142}{53}$, $v=\frac{142}{23}$,
$$ \{a,b,c,d\} = \left\{ \frac{19635}{6532}, -\frac{46592463}{201832268},
\frac{84196064}{50458067}, -\frac{1144273}{8775316} \right\} $$
and $t=\frac{59}{4}$, $v=\frac{59}{34}$,
$$ \{a,b,c,d\} = \left\{ \frac{2325}{4012}, \frac{187020623}{9949760}, \frac{261411943}{9949760},
\frac{13104399}{146320} \right\}. $$

\section{Torsion $\mathbb{Z}/2\mathbb{Z} \times \mathbb{Z}/4\mathbb{Z}$} \label{sec:z2z4}

Consider the curve
$$
E: \quad y^2 = x (x+p_1) (x+p_2), $$
where
$p_1=(b-a)(d-c)$, $p_2=(c-a)(d-b)$.

By the $2$-descent (see \cite[Theorem 4.2]{Knapp}), the point
$$ Q=((ad+1)(bc+1),\sqrt{(ab+1)(ac+1)(ad+1)(bc+1)(bd+1)(cd+1)}) $$
is in
$2E(\mathbb{Q})$, because
$x_1=(ad+1)(bc+1)$, $x_1+p_1=(ac+1)(bd+1)$ and $x_1+p_2=(ab+1)(cd+1)$ are perfect squares.
The point $Q$ is of order $2$ if $ad+1=0$, i.e. if $d=-1/a$.
So the point $R$ such that $2R=Q$ is of order $4$
and $E(\mathbb{Q})$ has a subgroup $\mathbb{Z}/2\mathbb{Z} \times \mathbb{Z}/4\mathbb{Z}$.

Thus we need to find rational Diophantine quadruples which contain a subtriple of the form
$\{a,-1/a,b\}$. But this is already done in \cite[Section 4]{D-glas}.
By the notation $\alpha=u, T=t$, we get
\begin{align*}
a &= \frac{ut+1}{t-u}, \\
b &= \Big(\Big(\frac{t+u}{t-u}\Big)^2-1\Big)/a = \frac{4tu}{(tu+1)(t-u)}.
\end{align*}
It remains to find the fourth element $c$ of the quadruple.
We may take $c$ such that $\{a,b,c,d\}$ is a regular quadruple,
$$ (a+b-c-d)^2 = 4(ab+1)(cd+1). $$
In that way, we get
$$ c= \frac{(u-1)(u+1)(t-1)(t+1)}{(ut+1)(t-u)}. $$

Hence, we showed that there are infinitely many rational Diophantine quadruples which induces
curves with torsion group $\mathbb{Z}/2\mathbb{Z} \times \mathbb{Z}/4\mathbb{Z}$.

Other possibility it is take $c$ to be
$\frac{8(d-a-b)(a+d-b)(b+d-a)}{(a^2+b^2+d^2-2ab-2ad-2bd)^2}$
(see \cite[Proposition 3]{D-quin}), so we get
{\footnotesize
$$ c=\frac{8(t-u)(ut+1)(-4ut+t^2+u^2+u^2t^2+1)(u-1)(u+1)(t-1)(t+1)(t^2+4ut+u^2+u^2t^2+1)}
{(1-8ut-12u^2t^2+2u^2+2t^2-8u^3t^3+u^4t^4+8ut^3+t^4+u^4+8tu^3+2u^2t^4+2u^4t^2)^2}. $$
}%

Within this two-parametric family of quadruples, we were able to find two examples of curves with rank
equal to 6,
for $(t,u)=(3,1/12)$ and $(t,u)=(25/2,31/8)$,
corresponding to the quadruples
$$ \left\{ \frac{32}{7}, \frac{15}{14}, \frac{4808055945}{27887662112}, -\frac{7}{32} \right\} $$
and
$$ \left\{ \frac{791}{138}, \frac{24800}{54579}, \frac{14188099227120}{9044268302161}, -\frac{138}{791} \right\}, $$
respectively.

\section{Torsion $\mathbb{Z}/2\mathbb{Z} \times \mathbb{Z}/6\mathbb{Z}$} \label{sec:z2z6}

To achieve the torsion group $\mathbb{Z}/2\mathbb{Z} \times \mathbb{Z}/6\mathbb{Z}$,
we will find quadruples for which the condition
$3Q=\mathcal{O}$, since then the point $R$ such that $2R=Q$ will be a point of order $6$.
To simplify the condition $3Q=\mathcal{O}$, we assume that $\{a,b,c,d\}$ is a regular quadruple.
Then we may take $R= \pm P$ by \cite[Proposition 1]{D-irr}.
There are several known parametrizations of Diophantine triples $\{a,b,c\}$ (see e.g. \cite{DKP-split,DP-racsam2,KNL}).
We will use the following parametrization due to Lasi\'c \cite{KNL}:
\begin{align*}
a &= \frac{2 t_1 (1 + t_1 t_2 (1 + t_2 t_3))} {(-1 + t_1 t_2 t_3) (1 + t_1 t_2 t_3)}, \\
b &= \frac{2 t_2 (1 + t_2 t_3 (1 + t_3 t_1))} {(-1 + t_1 t_2 t_3) (1 + t_1 t_2 t_3)}, \\
c &= \frac{2 t_3 (1 + t_3 t_1 (1 + t_1 t_2))} {(-1 + t_1 t_2 t_3) (1 + t_1 t_2 t_3)},
\end{align*}
followed by the substitutions which were already used in \cite{DP-racsam2}:
\begin{align*}
t_1 &= \frac{k}{t_2 t_3}, \\
t_2 &= m- \frac{1}{t_3}.
\end{align*}
From the regularity equation
$$ (a+b-c-d)^2 = 4(ab+1)(cd+1), $$
we get as one solution for $d$:
$$ d= \frac{-2(1-t_1+t_3t_1)(-t_3+t_2t_3+1)(-t_2+1+t_1t_2)(-1+t_1t_2t_3)}{(1+t_1t_2t_3)^3}. $$
Now the condition $3Q=\mathcal{O}$ become a complicated algebraic equation in terms of $t_3$, $k$ and $m$.
We searched for numerical solutions of this equation with $t_3$ and $k$ of small heights,
and tried to figure out some regularities within the solutions.
We have noted that there are solutions satisfying
$$ 2+t_3k = \frac{3}{k^2+2}, $$
i.e.
$$ t_3=\frac{-(2k^2+1)}{k(k^2+2)}. $$
By inserting this in $3Q=\mathcal{O}$, and factorizing the obtained expression,
we get a factor
$$ 4k^4m-4k^3m+6k^2m-2km+2m+6k-3k^2+6k^3-3k^4, $$
and thus we can express $m$ in terms of $k$ as
$$ m = \frac{3k(-2+k-2k^2+k^3)}{2(2k^2+1)(k^2-k+1)}. $$
Finally, we can express $a,b,c,d$ in terms of $k$:
\begin{align*}
a &= \frac{-2k(k^2+2)(3k^3-2k^2+2k-2)}{(1+k)(k-1)(2k^2+1)(2k^2+k+2)}, \\
b &= \frac{-k(1+k)(k-1)(2k^2+k+2)(4k^2-k+4)}{2(k^2+2)(k^2-k+1)^2(2k^2+1)}, \\
c &= \frac{2(2k^2+1)(2k^3-2k^2+2k-3)}{(1+k)(k-1)(2k^2+k+2)(k^2+2)}, \\
d &= \frac{(2k^2+1)(1+k)(k^2+2)(k-1)}{2(k^2-k+1)^2(2k^2+k+2)}.
\end{align*}
Hence, we obtained an infinite family of quadruples which induces
curves with torsion group $\mathbb{Z}/2\mathbb{Z} \times \mathbb{Z}/6\mathbb{Z}$.

We tried to compute the rank of curves corresponding to parameters $k$ with small numerators and denominators
by {\tt mwrank} \cite{Cremona} and {\tt magma} \cite{Magma},
and we found two curves with rank equal to $3$, for $k=23$ and $k=-\frac{22}{13}$,
corresponding to the quadruples
$$ \left\{ -\frac{16051953}{11214104}, -\frac{170244712}{1784519841}, \frac{914623}{5622936}, \frac{5498328}{10310521} \right\} $$
and
$$ \left\{ -\frac{18873668}{3382575}, \frac{821921100}{5086844387}, -\frac{26226421}{4890900}, \frac{1090383}{6661892} \right\}, $$
respectively.

\section{Torsion $\mathbb{Z}/2\mathbb{Z} \times \mathbb{Z}/8\mathbb{Z}$} \label{sec:z2z8}

Let us start with the quadruple $\{a,b,c,d\}$ from Section \ref{sec:z2z4}:
\begin{align*}
a &= \frac{ut+1}{t-u}, \\
b &= \frac{4ut}{(ut+1)(t-u)}, \\
c &= \frac{(u-1)(u+1)(t-1)(t+1)}{(ut+1)(t-u)}, \\
d &= -\frac{t-u}{ut+1}.
\end{align*}
The point $Q=(0,0)$ of order $2$ is in $2E(\mathbb{Q})$, so there
is point $R$ such that $2R=Q$.
To get the torsion group $\mathbb{Z}/2\mathbb{Z} \times \mathbb{Z}/8\mathbb{Z}$,
we just have to force the point $R$ to be in $2E(\mathbb{Q})$, i.e. $R=2S$,
 for a point $S$ in $E(\mathbb{Q})$. Then the point $S$ will be of order $8$.
The point $R$ has coordinates:
$$ \left( \frac{(u+t)^2 (u t-1)^2}{(ut+1)^2 (t-u)^2},
\frac{(u^2+1) (t^2+1) (u+t)^2 (u t-1)^2}{(ut+1)^3 (t-u)^3} \right). $$
We want that for $x_1=\frac{(u+t)^2 (u t-1)^2}{(ut+1)^2 (t-u)^2}$,
$x_1$, $x_1+p_1$ and $x_1+p_2$ are squares.
But $x_1$ is already a square, while $x_1+p_1=\Box$ and $x_1+p_2=\Box$ both lead to the
same condition:
$(u^2+1)(t^2+1)$ is a square.
Put
$$ (u^2+1)(t^2+1) = (u^2+1+(t-u)v)^2 $$
and we get
$$ t = -\frac{-2u^2v-2v+v^2u+u^3+u}{u^2+1-v^2}, $$
and finally,
{\footnotesize
\begin{align*}
a &= \frac{(u-v+1)(u-v-1)}{2(u-v)}, \\
b &= -\frac{2(u^2+1-v^2)u(-2u^2v-2v+v^2u+u^3+u)}{(u^2+1)^2(u-v)(u-v+1)(u-v-1)}, \\
c &= \frac{(-2u^2v-2v+v^2u+u^3+u-u^2-1+v^2)(-2u^2v-2v+v^2u+u^3+u+u^2+1-v^2)(u+1)(u-1)}{2(u^2+1)^2(u-v)(u-v+1)(u-v-1)}, \\
d &= -\frac{2(u-v)}{(u-v+1)(u-v-1)}.
\end{align*}
}%
Here $u$, $v$ are arbitrary rationals such that
{\footnotesize
\begin{align*}
& uv(u+1)(u-1)(-uv+v+1+u^2)(-uv+1+u^2)(-uv-v+1+u^2)(u-v+1)(u-v)(u-v-1) \\
& \times (u^2+1-v^2) (-2u^2v-2v+v^2u+u^3+u)(-2u^2v-2v+v^2u+u^3+u-u^2-1+v^2) \\
& \times (-2u^2v-2v+v^2u+u^3+u+u^2+1-v^2)
\end{align*}
}%
 is non-zero (so that $a,b,c,d$ are distinct non-zero rationals).

\medskip

By comparing $j$-invariant of our curve $E$
with the general form of an elliptic curve with torsion group
$\mathbb{Z}/2\mathbb{Z} \times \mathbb{Z}/8\mathbb{Z}$, i.e.
$$ y^2 = x \Big(x+ \Big(\frac{2T}{T^2-1}\Big)^2\Big) \Big(x+ \Big(\frac{T^2-1}{2T}\Big)^2\Big), $$
we find that $j$-invariants coincide for $T=\frac{v}{vu-u^2-1}$.
This implies that every curve with torsion group
$\mathbb{Z}/2\mathbb{Z} \times \mathbb{Z}/8\mathbb{Z}$ can be obtained
from a rational Diophantine quadruple.

For example, the quadruple
$$ \left\{ \frac{1804}{1197}, -\frac{226796}{539847}, \frac{303199}{239932}, -\frac{1197}{1804} \right\} $$
gives the curve with torsion $\mathbb{Z}/2\mathbb{Z} \times \mathbb{Z}/8\mathbb{Z}$ and rank $3$
(equivalent to the curve found by Connell and Dujella in 2000, see \cite{D-torsion}),
what is the largest known rank for curves with this torsion group.

\bigskip

{\bf Acknowledgement.} The authors would like to thank Matija Kazalicki for useful
comments on the previous version of this paper.

\end{document}